
\documentstyle{amsppt}
\magnification 1200
\def\boxit#1{\vbox{\hrule\hbox{\vrule\kern2.5pt\vbox{\kern2.5pt#1
\kern2.5pt}\kern2.5pt\vrule}\hrule}}
\topmatter
\title A few general points of the Hilbert scheme of K3 surfaces \endtitle
\author Francisco Javier Gallego \\ and \\ B. P. Purnaprajna
\endauthor
\affil Brandeis University \\ Department of Mathematics \endaffil
\address{Francsico Javier Gallego: Departamento de Algebra,
 Facultad de Matematicas,
 Universidad Complutense de Madrid, 28040 Madrid,Spain}\endaddress
\email{gallego\@sunal1.mat.ucm.es}\endemail
\address{ B.P.Purnaprajna:
Dept of Mathematics
 Brandeis University,
Waltham MA 02254-9110}\endaddress
\email{purna\@max.math.brandeis.edu}\endemail
\thanks{We are happy to
thank our advisor David Eisenbud for his help, patience
encouragement. We would also like to thank Andrea Bruno and Enrique Arrondo
for helpful discussions}\endthanks
\endtopmatter
\document
\vskip .3 cm

\roster
\item ""{\bf  Introduction}
\item "{\bf 1.}" {\bf K3 carpets}
\item "{\bf 2.}" {\bf A characterization for ribbons}
\item "{\bf 3.}" {\bf Smoothings of K3 carpets}
\item "{\bf 4.}" {\bf The Hilbert scheme
near the point of a K3 carpet }
\endroster

\heading  Introduction\endheading

This article deals with the study of K3
carpets.
D. Bayer and D. Eisenbud say in [BE] that ``{\it a ribbon}''
(supported on
$\bold P ^1$ inside $\bold P^{g-1}$ and with arithmetic genus $g$)
``is the answer to the riddle: What is the limit of the canonical
model of a smooth curve as the curve degenerates to a hyperelliptic
curve?'' Analogously one would say that a {\it K$3$ carpet} is the
answer to the riddle: What is the limit of the embedded model of
a smooth polarized K3 surface as the polarized surface
degenerates to a hyperelliptic polarized surface? To justify this
claim we devote much of this article.

K3 carpets possess some interesting features. On the one hand there are
few of them. In Section 1 we see that there is only one K3
carpet supported on a given rational normal scroll (in the same
way as a canonical ribbon is a double structure on a rational normal
curve, the reduced structure of a K3 carpet is a rational normal
scroll). Thus one can in some sense think of the set of all K3
carpets as something discrete. On the other hand, some of them,
(the ones supported on ``balanced" scrolls) are still general, in
the sense that they are smooth points of the Hilbert scheme. Hence
K3 carpets form a small class of very degenerate objects (they are
nowhere reduced and one step more degenerate than such reduced
non-normal K3 surfaces as the unions of two (distinct) rational
normal scrolls) which are nevertheless general.

Another interesting feature is that the hyperplane section of a K3
carpet is a canonical ribbon. The
study of canonical ribbons has been proposed by Bayer, Eisenbud,
Green and Schreyer among others as a means to solve the so-called
Green's conjecture. Briefly, in its original form this conjecture
relates the graded Betti numbers of the minimal free resolution of a
canonical curve to the Clifford
index of the curve (the Clifford index of a smooth curve of genus
$g \geq 3$ is defined as the minimum, over all line bundles $L$ on
the curve such that $\text h^0(L) >1$ and $\text h^1(L) >1$ of the
quantity $\text{Cliff}L=\text{deg}L-2(\text h^0(L)-1)$. More
loosely said, the Clifford index tells us how special the most special
line bundle which the curve possesses is). More precisely one
expects that  the canonical bundle will satisfy the
 property $N_{p}$ but not the  property $N_{p+1}$ iff $p$ is the
Clifford index of the curve. Thus Green's conjecture generalizes
classical results by Noether and Petri
(c.f. [ACGH]; for details on Green's conjecture see [E]). Results by
Eisenbud and Green [EG] and Fong [F] yield that an affirmative
answer to Green's conjecture in the case of canonical ribbons will
imply Green's conjecture for general curves. Since K3 carpets are
arithmetically Cohen-Macaulay, the Betti numbers of the minimal free
resolution of a K3 carpet are the same as the Betti number of the
hyperplane section, a canonical ribbon. Progress toward the
computation of the minimal resolution of a K3 carpet has been
made by Dave Bayer and David Eisenbud in [BE] where they
compute the graded Betti numbers of a nonminimal resolution of a K3
carpet.

Our work on K3 carpets focuses on answering two questions: Are these
objects  smoothable ? and Do they correspond to smooth points of the
Hilbert scheme?. The first question is dealt with
and answered affirmatively in Section 3.
To prove this result we use the idea, already introduced,  that
a (suitable) K3 carpet is morally the ``image" of the morphism
associated to a hyperelliptic linear system. To show this we use a
characterization given in Section 2 which allows us to decide by
induction on the dimension (cutting with a hyperplane) whether a
scheme is a ribbon. We also use properties of hyperelliptic linear
systems on K3 surfaces and of the moduli of K3 surfaces.
We would
like to point out here one difference between the case of canonical
ribbons and the case of K3 carpets. While canonical ribbons can be
thought as ``canonical models" of hyperelliptic curves, not all K3
carpets are ``models" of smooth hyperelliptic K3 surfaces.
More precisely, rational normal scrolls with a rational curve
with low self-intersection cannot be realized as images of
morphisms associated to hyperelliptic linear systems, and
hence the riddle posed before does not make sense for them.
However we are able to prove that this kind of K3 carpets are
smoothable by showing that they deform to more general ones.

In Section 4 we deal with the study of the Hilbert scheme near the
locus of K3 carpets. Our main result is that K3 carpets supported
on ``balanced" scrolls are smooth points of the Hilbert scheme. Here
another departure from the case of ribbons occurs. While both
K3 carpets and canonical ribbons are smoothable (i.e, both belong
to the component parametrizing smooth varieties in their
respective Hilbert schemes), contrary to the case of canonical
ribbons, not all K3 carpets are smooth points of the Hilbert
scheme (some of them even belong to several components of the
Hilbert scheme as noted in Theorem 4.3).

\heading {1. K3 carpets}\endheading
\proclaim{Conventions} \rm {Throughout this article we work over
$\bold C$. A rational normal scroll or simply a scroll will
always mean a smooth rational normal scroll of dimension $2$. We
will denote by $F_n$ the rational ruled surface whose minimal
section has self-intersection $-n$.}
\endproclaim
In this section we introduce our main objects of study, the K3
carpets and some properties of them, which we will use later in
the article. We start with a couple of definitions:

\proclaim {Definition 1.1 ([BE], \S 1)}
A double structure or a
ribbon on a reduced connected scheme $D'$ is a scheme $D$
equipped with an isomorphism $D' \to D_{red}$, such that the
ideal sheaf $\Cal I$ of $D'$ in $D$ satisfies $\Cal I^2 =0$ and
is a line bundle on $D'$.
\endproclaim

\proclaim {Definition  1.2} A
K3 carpet $\tilde S$ is a double
structure on a rational normal
scroll $S$ (i.e., a double structure
embedded in some
$\bold P^g$, whose reduced
structure is a rational normal
scroll in $\bold P^g$) such that
its dualizing sheaf
$\omega_{\tilde S}$ is trivial and
$h^1({\Cal O}_{\tilde S})=0$.\endproclaim

\vskip .5 cm
An important
fact about K3 carpets (which will be certainly instrumental to
our proof of the main result of this article, namely, the
smoothing of K3 carpets) is stated in this
\proclaim{Theorem 1.3} There is a unique K3 carpet (up to
multiplication by  scalar)
on a given rational normal scroll.
\endproclaim

Before we prove Theorem 1.3 we need to state two
lemmas  which  are variants of results in [HV]. The lemmas identify
the conormal bundle of the reduced structure of the K3 carpet inside
the carpet itself. From them it follows that the K3 double
structures on a scroll $S$ in
$\bold P^g$ corresponds to the global sections of a twist of the
normal bundle of $S$. The proofs use the same ideas of [HV].
\proclaim{Lemma 1.4} Defining a double structure $\tilde S$
on a smooth subvariety $S$ of a smooth variety
$Z$ is equivalent to
giving a line subbundle $\Cal L$
of
$\Cal N_{S/{Z}}$. This line
bundle $\Cal L$ is the normal
bundle of $S$ in $\tilde S$ \endproclaim
{\it Proof.} Let ${\Cal L}
\subseteq \Cal N_{S/{Z}}$ be a line bundle and
${\Cal I}={\Cal I}_{Z}(S)$ the ideal sheaf of $S$
in $Z$. Let
$w$ be the  surjective homomorphism:$$  w:
{\Cal I}  \rightarrow {\Cal I}/{\Cal I}^2=
\Cal N_{S/Z}^* \rightarrow {  {\Cal L
}}^*$$ Let ${\Cal J}=\ker w$. The ideal sheaf  ${\Cal J}$
defines a subscheme $\tilde S$ in $Z$. From the
exact sequences $$0 \to {\Cal J} \to {\Cal I} \rightarrow {
{\Cal L
}}^* \rightarrow 0$$ and
$$\matrix
0 & \to & \Cal I/ \Cal J & \to &\Cal O _{Z} /\Cal J & \to & \Cal O _{Z} /
\Cal I & \to & 0 \\
  &     &  \|            &     &   \|                       &     &  \|
                 &     &   \\
0 & \to & {\Cal L}^* & \to &\Cal O _{\tilde S}          & \to & \Cal O
_S                    & \to & 0
\endmatrix$$
we see that the ideal sheaf defining $S$ in $\tilde S$ is the
line bundle ${ {\Cal L }}^*$. Since ~{$({ {\Cal L
}}^*)^{ 2}=({\Cal
I}/{\Cal J})^{ 2}= 0$} by
construction of ${\Cal J}$, it follows that
$\tilde S$ is a ribbon.

Conversely, let $\tilde S$ be a
double structure embedded in
$Z$, let $S$ be its reduced
part, and let ${\Cal I}$,${\Cal J}$ be
their respective ideal sheaves
in $Z$. By the definition of
ribbon ${\Cal I}^{ 2}
\subseteq {\Cal J}$, so
\text{$\Cal N_{S/{Z}}={\Cal I}/{\Cal
I}^2$} surjects onto ${\Cal
I}/{\Cal J}$, which is the
conormal bundle of $S$ in
$\tilde S$, in particular, a
line bundle. \boxit{}

\vskip .5 cm

\proclaim{Lemma 1.5} Let $S$
be a rational normal scroll and
$\tilde S$ a carpet whose
reduced part is $S$ and let
$\Cal L$ be the dual of the
ideal sheaf defining $S$ in
$\tilde S$. Then $\tilde S$ is a
K3 carpet if and only if $\Cal L
\simeq { {{\omega}_S}}^*$. \endproclaim

{\it Proof}.
First assume ${\Cal L}
\simeq \omega_S^*$, so we have an exact sequence
$$ 0 \rightarrow
{\omega}_S \rightarrow {\Cal
O}_{\tilde S} \rightarrow {\Cal
O}_S \rightarrow 0  \eqno (1.5.1) $$
{}From $\text H^1({\Cal
O}_S)=\text H^1({\omega}_{
S})=0$ it follows that
$\text H^1({\Cal
O}_{\tilde S})=0$. If we apply
to (1.5.1) the functor ${\Cal
Ext} _{{\Cal O}_{\bold
P^g}}^{g-2}(-,{\omega}_{\bold P^g})$
 we get (see, for example, [H], p. 235)
 $$0 \rightarrow
{\omega}_S \rightarrow
{\omega}_{\tilde S} \rightarrow
{\Cal O}_S \rightarrow 0 .$$
Since
$\text H^1 ({\omega}_S)=0$, the map $\text H^0
({\omega}_{\tilde S}) \rightarrow
\text H^0({\Cal
O}_S)
$ is an epimorphism.
Therefore $1 \in H^0({\Cal
O}_S)$ can be lifted to $H^0
({\omega}_{\tilde S})$ and hence
${\omega}_{\tilde S} \simeq {\Cal
O}_{\tilde S}=0$.

Now assume $\tilde S$ is a K3
carpet. Apply the functor ${\Cal
Ext}_{{\Cal O}_{\bold P^g}}^{g-2}(-,
{\omega}_{\bold P^g})$
to the exact sequence $$ 0
\rightarrow \Cal I \rightarrow
{\Cal O}_{\tilde S} \rightarrow {\Cal
O}_S \rightarrow 0$$ to obtain
$$0
\rightarrow {\omega}_S \rightarrow
{\omega}_{\tilde S} \rightarrow
 {\Cal I ^*}
 \otimes {\omega}_S \rightarrow
0 \eqno (1.5.2)$$
If we tensor
(1.5.2) with ${\Cal
O}_S$ and we use the fact that ${\omega}_{\tilde S} \simeq {\Cal
O}_{\tilde S}$ we get a surjection ${\Cal
O}_S \rightarrow  {\Cal I ^*}
\otimes {\omega}_S $. Thus   $
{\Cal I ^*}
\otimes {\omega}_S \simeq {\Cal
O}_S$ and $\Cal I \simeq
{\omega}_S$. \boxit{}

\vskip .5 cm

(1.6) Lemmas 1.4 and 1.5 imply that in order to
see how many K3 carpets are
supported in a particular
rational normal scroll $S$, one
has to compute how many bundle
inclusions there are from
$ \omega_S^*$ into
$\Cal N_{S,{\bold P^g}}$, or equivalently, how many nowhere
vanishing sections there are
in $\text H^0 (\Cal N_{S,{\bold P^g}} \otimes
{\omega}_S)$. Therefore Theorem 1.3 follows from the following

\vskip .6 cm

\proclaim{Proposition 1.7} Let $S = S(a,b)$, $a \geq b$ be the rational normal
scroll that corresponds to the embedding of ${\bold P}({\Cal E})$ into ${\bold
P}^N$ by ${\Cal O}_{{\bold P}({\Cal E})} (1)$, where  ${\Cal E} = {\Cal
O}_{{\bold P}^1}(a) \oplus {\Cal O}_{{\bold P}^1}(b)$. Let $\omega$ be the
canonical bundle of ${\bold P}({\Cal E})$. Then $\text H^0 (\Cal N_{S/{{\bold
P}^N}}
\otimes \omega) = \bold C <s>$, where s is a nowhere vanishing section, and
\text{$\text H^1 (\Cal N_{S/{{\bold
P}^N}}
\otimes \omega) = \text H^2 (\Cal N_{S/{{\bold
P}^N}}
\otimes \omega) = 0$.}
\endproclaim
{\it Proof}. We use the exact sequence
$$ 0 \to {\Cal T}_S \otimes \omega \to {\Cal T}_
{{\bold P}^N}|_S \otimes \omega \to
\Cal N_{S/{\bold P}^N} \otimes \omega \to 0
\eqno (1.7.1)$$
to compute $\pi _* (\Cal N_{S/{\bold P}^N} \otimes
\omega)$, where $\pi$ denotes the projection from $S$ to $\bold P
^1$. To compute $\pi _* ({\Cal T}_S \otimes \omega)$  we use
the exact sequences
$$ 0 \to \Cal T _{S /_{{\bold P}^1}} \otimes \omega \to {\Cal T}_S
\otimes \omega \to \pi ^* {\Cal T}_{{\bold P}^1} \otimes \omega
\to 0 $$
and $$ 0 \to \Omega _{S /_{\bold
P^1}} \to \pi ^*
 \Cal E \otimes \Cal O _{{\bold P}(\Cal E )} (-1)
 \to
 \Cal O _{{\bold P}(\Cal E )} \to 0   \eqno (1.7.2)
 $$
which is a relative version of the  Euler
sequence.

Let $\Cal E ' = \Cal E \otimes \Cal O _{\bold P ^1} (-a) = \Cal O
_{\bold P ^1}
\oplus \Cal O _{\bold P ^1} (b-a)$. Then
$\Cal O _{{\bold P}(\Cal E ' )} (1) = \Cal O
 (C_0)$, where $C_0$ denotes the
 minimal section of $\pi :S \to
\bold P ^1$ . Therefore, by exact sequence (1.7.2)
we obtain
$$\Omega _{S {/_{\bold P ^1}}} =\bigwedge ^2 (\Cal O
(-C_0) \oplus \Cal O  (-C_0+(b-a)f))=\Cal O
(-2C_0+(b-a)f)$$ Also, we know that $\omega =
\Cal O  (-2C_0 + (b-a-2)f)$ and that
$\pi ^*  \Cal T_{\bold P}^1 = \Cal O
 (2f)$. Hence we obtain the sequence
$$ 0 \to \Cal O (-2f) \to  {\Cal T}_S \otimes \omega \to
\Cal O (-2C_0 +(b-a)f) \to 0 \ .$$
We apply $\pi _*$  and get:
$$ \displaylines{0 \to {\Cal O}_{{\bold P}^1} (-2) \to \pi _*
({\Cal T}_S \otimes \omega) \to \pi _* (\Cal O
(-2C_0 +(b-a)f) =0 \to \cr \to R^1 \pi _*
{\Cal O} (-2f ) =0 \to R^1 \pi _* ({\Cal T}_S
\otimes \omega) \to R^1 \pi _* (\Cal O (-2C_0 +(b-a)f) \to 0 \ .}$$
Therefore $\pi _* ({\Cal T}_S \otimes \omega) = {\Cal O}_{{\bold P}^1} (-2)$
and
$$R^1 \pi _* ({\Cal T}_S \otimes \omega) = R^1 \pi _* (\Cal O (-2C_0 +(b-a)f)
= (\pi _* \Cal O) ^* = \Cal O _{\bold P ^1}$$ by relative Serre duality.

To compute $\pi _* ({\Cal T}_{{\bold P}^N}|_S
\otimes \omega)$  we push forward
the presentation of ${\Cal T}_{{\bold P}^N}|_S
\otimes \omega$, which comes from the Euler
sequence:
$$ 0 \to \Cal O (-2C_0 +(b-a-2)f) \to (\Cal O (-C_0
+(b-2)f)) ^{\oplus (N+1)} \to
{\Cal T}_{{\bold P}^N}|_S \otimes \omega \to 0$$
and we obtain:
$$ \displaylines{0 \to \pi _* (-2C_0 +(b-a-2)f)
\to \pi _*
[(\Cal O (-C_0 +(b-2)f)) ^{\oplus (N+1)}
] =0 \to \cr
\to \pi _* ({\Cal T}_{{\bold P}^N}|_S \otimes
\omega)
\to
  R^1 \pi _* (-2C_0 +(b-a-2)f) \to \cr
  \to R^1 \pi _* [(\Cal O
(-C_0 +(b-2)f)) ^{\oplus (N+2)} ] =0} $$

 Thus $R^1 \pi _* ({\Cal T}_{{\bold P}^N}|
_S \otimes \omega) =0 $ and $$ \pi _* ({\Cal T}_
{{\bold P}^N}|_S
\otimes \omega) = R^1 \pi _*
(-2C_0 +(b-a-2)f) =
 (\pi _* \Cal O (2f))^* = \Cal O
_{{\bold P} ^1} (-2) \ ,$$   by relative
Serre duality.
Applying $\pi _*$ to  (1.7.1) we get:
$$ 0 \to \Cal O _{\bold P ^1} (-2) \to \Cal O _
{\bold P ^1} (-2) \to \pi _* (\Cal N_{S/{\bold P}^N} \otimes
\omega)
\to \Cal O
_{\bold P ^1} \to 0 $$
Hence
$$\displaylines{\rlap{(1.7.3)}\hfill \pi _* (\Cal N_{S/{\bold P}^N} \otimes
\omega)
 = \Cal O _{\bold P ^1} \hfill\llap{\text{and}\qquad}\cr
\rlap{(1.7.4)}\hfill R^1\pi _* (\Cal N_{S/{\bold P}^N} \otimes
\omega)
 = 0 \ . \hfill\cr} $$
This means that there exists a nonzero global
section $s$ of $\Cal N_{S/{\bold P}^N} \otimes \omega$. This section cannot
vanish identically at any fiber of $\pi$. But the fibers of $\pi$ are
projective lines and hence, by (1.7.3), the restriction of $\Cal N_{S/{\bold
P}^N}
\otimes
\omega$ to a fiber is isomorphic to $\Cal O_{\bold P ^1} \oplus
\Cal F$ for
some vector bundle $\Cal F$ without global sections. This implies
that the restriction of
$s$ to each fiber is nowhere vanishing. This proves the statement about $\text
H^0(\Cal N_{S/{\bold P}^N} \otimes \omega)$. The statement about $\text H^1$
and
$\text H^2$ follows from (1.7.3) and (1.7.4).
\boxit{}

\heading 2. A characterization for ribbons \endheading

The next theorem gives a way to decide whether a scheme is a ribbon
by using induction on the dimension.

\proclaim {Theorem 2.1}Let $D$ be an  scheme such that
$D_{red}$ is equidimensional. $D$
is  a ribbon iff for every closed point $p
\in  D_{red}$ there exists  $h_p \in \Cal O_{D,p}$ such
that $\Cal O_{D,p} /(h_p)$ is the structure sheaf of
a ribbon of dimension dim$D - 1$, whose reduced structure is  $\Cal
O_{D_{red},p} /(h_p)$.\endproclaim


In order to  prove the theorem we will need the following

\proclaim {Remark 2.1.1} Let $M$ be a module over a ring
 $A$. Let $a \in A$ be a non-zero-divisor for $M$. Then
Tor$^1_A (M,A/(a)) = 0$. For example, in the
situation of Theorem 2.1, $h_p$ is non-zero-divisor for
$\Cal O _{D_{red}}$ and $ \text{Tor} ^1 (\Cal O _{D_{red}} , \Cal O_{
D,p} / (h_p)) = 0$\endproclaim

(2.2) {\it Proof of the Theorem 2.1.}
The ``only if" part is trivial.
For the ``if"
part  let $\Cal I$ be the ideal sheaf of $D_{red}$ in $D$.
We have to show that   $\Cal I / \Cal I ^2$
is a locally free sheaf over $\Cal O _{D_{red}}$
and that $\Cal I ^2 =
0$.

\vskip .25 cm
{\it Step 1} ($\Cal I / \Cal I ^2$
 is locally free).

Let us fix a closed point $p \in D_{red}$.  The ideal $(\Cal I
 _p + h_p)/(h_p) = \Cal J _p$ is
the ideal of the reduced part of a ribbon,  so $\Cal J _p/\Cal J _p
^2 = \Cal J_p$ is a free module generated by one
element. By Remark
 2.1.1, $\Cal J _p = \Cal I_p \otimes
\Cal O_{D,p} /(h_p)$ and hence, $\Cal I _p/\Cal I ^2
_p \otimes \Cal O_{D_{red},p} /(h_p) = \Cal I _p/
\Cal I ^2
_p \otimes \Cal O_{D,p} /(h_p) = \Cal J _p/\Cal J _p
^2$ is a free cyclic module. Thus, by Nakayama, $\Cal
I_p
/\Cal I_p ^2$ is also a cyclic module over $\Cal O
_{D_{red},p}$ generated by an element not vanishing
at $p$. This is true for any closed point
 $p \in D_{red}$, i.e.,
the rank of $\Cal I/\Cal I ^2$ is 1 for any closed
point $p \in D_{red}$. Hence $\Cal I/\Cal I ^2$ is
locally free over $\Cal O_{D_{red}}$.
\vskip .25 cm
{\it Step 2}
$(\Cal I ^2 = 0)$.

Fix again $p \in D_{red}$. We claim that $\Cal I _p
^2 \otimes \Cal O _{D,p}/(h_p) = 0$.
Indeed. As we remark in Step 1, $(\Cal I
 _p + h_p)/(h_p) = \Cal J _p$ is
the ideal of the reduced part of a ribbon in its
structure sheaf, so $\Cal J_p^2 = 0$ and hence $\Cal
I_p^2 \subset (h_p)$. This implies that
$\Cal O _{D,p}/\Cal I_p^2 \otimes \Cal O_{D,p}/(h_p)
= \Cal O_{D,p}/(h_p)$, so tensoring

$$0 \to \Cal I_p ^2 \to \Cal O_{D,p} \to \Cal O_{D,p} /\Cal I_p ^2 \to 0 $$
by $\Cal O_{D,p}$ we get:
$$\displaylines{
0 \to \text {Tor} ^1 (\Cal O_{D,p} /\Cal I_p ^2 ,
\Cal O_{D,p} /(h_p)) \to \Cal I_p ^2
\otimes \Cal O_{D,p}/(h_p) \cr
\to \Cal O_{D,p} /(h_p) @>\simeq >>
\Cal O_{D,p} /(h_p) \to 0 \cr}  $$
In order to prove our claim, it suffices
to prove that $\text {Tor} ^1 (\Cal O_{D,p} /\Cal I_p ^2 ,
\Cal O_{D,p} /(h_p)) = 0$.  If we tensor
$$  O \to \Cal I _p /\Cal I_p ^2 \to \Cal O _{D,p}/ \Cal I_p ^2
\to \Cal O _{D_{red},p} \to 0 $$
by $\Cal O_{D,p} /(h_p)$ , we
obtain:
$$\displaylines{\text {Tor} ^1 (\Cal I_p /\Cal I_p ^2 ,
\Cal O_{
D,p} / (h_p)) \to \text{Tor} ^1 (\Cal O_{D,p}
 /\Cal
I_p ^2 , \Cal O_{
D, p} / (h_p))\cr
 \to
 \text{Tor} ^1 (\Cal O _{D_{red}} ,
\Cal O_{
D, p} / (h_p))\ .}$$

By Step 1, we know that $ \Cal I_p /\Cal I_p ^2 = \Cal O _{D_
{red},p}$, so again by Remark 2.1.1, \linebreak
 both   $ \text{Tor} ^1 (\Cal O
_{D_{red}} , \Cal O_{
D,p} / (h_p))$ and $ \text {Tor} ^1 (\Cal I_p /\Cal I_p ^2 ,
 \Cal O_{
D,p} / (h_p))$ vanish.
Therefore $\Cal I _p
^2 \otimes \Cal O _{D,p}/(h_p) = 0$ and by Nakayama's
lemma, $\Cal I_p ^2
 = 0$.  \boxit{}

\heading {3. Smoothings of K3 carpets} \endheading
The purpose of this section is to prove the existence of smoothings
of K3 carpets. By smoothing we mean a flat family over a smooth
curve with smooth generic fiber and with a special closed fiber
isomorphic to the K3 carpet. We prove the result in two
steps. Using the fact that rational normal scrolls  $F_0,
\dots , F_4$ admit a generically $2:1$ map from a
hyperelliptic K3 surface, we construct, in a rather explicit
way,  smoothings  of K3 carpets supported on  $F_0,
\dots , F_4$. Then, in Theorem 3.6 we will see that the remaining
K3 carpets lie in the closure of the locus parametrizing K3 carpets
supported on $F_0,
\dots , F_4$. In order to prove these results we will need some
auxiliary lemmas.

In this section, a smooth curve will mean either an algebraic
smooth curve or the analytic disc $\Delta$.

\proclaim {Remark 3.1} Let $\Cal X$ be a flat family over a smooth
 curve $T$. If
$\phi : \Cal X
\to \Cal Z$ is a  morphism  over $T$, then the image $\Cal Y$ of
$\phi$,
 is flat over $T$.\endproclaim

{\it Proof.}
Let $\pi$ be the morphism from $\Cal Y$ to $T$. By assumption $\pi _*
\phi _*
\Cal O _\Cal X$ is flat over $\Cal O_\Cal Y$ and therefore
\text{$\pi _* \Cal O _\Cal Y \hookrightarrow \pi _* \phi _*
\Cal O _\Cal X$} is a subsheaf of a torsion free sheaf on
$T$, so it is itself torsion free and hence flat. \boxit{}

\proclaim {Lemma 3.2} Let $\Cal X$ be a flat family of irreducible
varieties over
a smooth curve $T$.
Let
$\zeta $ be  a relatively globally generated, invertible
sheaf on $\Cal X$. Let $\phi $ be the morphism from $\Cal X$
to $\bold P_T^n$ induced by its relative complete
linear series and let $\Cal Y$ be the image of $\Cal X$ by $\phi$.
Assume that $\phi$ is an embedding outside
the central fiber, and a finite morphism of degree $2$
when restricted to the central fiber.
 Let $H$ be an hyperplane in $\bold P^n$. Then $(\Cal Y
\cap (H \times T))_t = \Cal Y_t \cap (H \times \{ t
\} )$ and $\Cal Y \cap (H \times T )$ is  flat
over $T$.\endproclaim

{\it Proof.}
The pull back of $H \times T$ is a Cartier divisor on $\Cal X$
whose zero locus defines a flat family of divisors $\Cal X'$.
Indeed, the only thing to be checked (c.f. [H], III.9.8.5)
is whether
the pullback of $H \times T$ is defined by a non-zero-divisor at
$\Cal O_ {\Cal X_t ,p}$, for all $t \in T$ and for all $p \in \Cal X_t$.
This is obvious, since $\Cal X_t$ is reduced and irreducible,
and $H$ does not contain $\phi _t (\Cal X_t)$.

Now, the image of $\Cal X'$ by $\phi$ is a flat family by the
previous observation. Hence, if we see that $\phi (\Cal X') =
\Cal Y \cap (H \times T )$, we are done. We have to prove that
the morphism $ \Cal O _{\Cal Y \cap (H \times T)}
\rightarrow  \phi_{0 _*}
\Cal O _\Cal X$, obtained by tensoring $\Cal O _\Cal Y \hookrightarrow \phi _*
\Cal O _\Cal X$ by $\Cal O_{\bold P_T^n}/\Cal I(H \times T)$,
is injective. Consider the exact sequence
$$ 0 \to  \Cal O _\Cal Y @> \alpha >> \phi _*
\Cal O _\Cal X \to \Cal F \to 0 \ .$$
 The rank of $\phi _*
\Cal O _\Cal X$ is 1 outside $\phi(\Cal X_0)$ and 2 at
 $\phi (\Cal X_0)$.  The injection $\alpha$ of $\Cal O _\Cal Y$ into
$\phi _*
\Cal O _\Cal X$ is given by a nowhere vanishing global section of
$\phi _*
\Cal O _\Cal X$; hence $\alpha$ is  an injection at each fiber. From
all  this, it follows that $\Cal F$ is supported at $\Cal Y_0$ and
has rank 1 at every closed point $y \in \phi (\Cal X_0)$ (i.e., it
is a line bundle on
 $\phi (\Cal X_0)$).
 By hypothesis $\phi (\Cal X_0)$ is an irreducible variety, so $H$
is locally a nonzero divisor at every point of $\phi (\Cal X_0)$.
Remark 2.2.1 implies that Tor$^1(\Cal F, \Cal O_{\bold P_T^n}/\Cal
I(H \times T)) = 0$, so
$ \Cal O _{\Cal Y \cap (H \times T)} \rightarrow
\phi_{0 _*}
\Cal O _\Cal X$ is injective as required.

The fact that $(\Cal Y
\cap (H \times T))_t = \Cal Y_t \cap (H \times \{ t
\} )$ is obvious.
\vskip .5 cm

We recover as a corollary of  Theorem 2.1, the following
result of Fong:

\proclaim{Corollary 3.3. ([F], Theorem 1.(i))} Let $\Cal C$ be a
flat family of smooth curves over a smooth curve $T$ such that its
central fiber is a hyperelliptic curve and its generic fiber is a
nonhyperelliptic curve. If $\Cal D$ is the image in $\bold
P^{2g-2}_T$ of $\Cal C$ by the relative complete linear series
of $\omega _{\Cal C /T}$, then the central fiber of $\Cal D$ is a
canonical ribbon.
\endproclaim

{\it Proof.}
Note first that by Remark 3.1 $\Cal
D$ is flat over $T$. Let $D$ be the central fiber of $\Cal
D$. We want to prove that $D$ is a canonical ribbon (recall that
$D$ is not the image of the central fiber of $\Cal C$).  The  degree
of
$D$ is
$2g-2$ and its arithmetic genus is
$g$. The reduced part of $D$, $D_{red}$, is a rational normal
curve. By Theorem 2.1, in order to see that $D$ is a ribbon we
need to check that at every point p of $D_{red}$, we can choose a
hyperplane
$H_p$  passing through $p$ such that $H_p \cap D$ is
isomorphic to $(g -1)$ copies of
 Spec$(\bold C (\epsilon ))$. To see this  choose $H_p$ through $p$
intersecting $D_{red}$ at $g-1$ distinct points. Lemma 3.2 tells
that
$H_p
\cap D$ is the flat limit of a family of $(2g-2)$ points. $H_p
\cap D$ must be non-reduced everywhere. If not, a point of $H_p
\cap D$ would be
a smooth point of $D$, and the degree of $D$ would be equal to the
degree of
$D_{red}$, which is
$g-1$. On the other hand, the degree of each component of
$H_p \cap D$ must be less or equal than two; otherwise there
would be reduced point in $H_p
\cap D$.
\boxit{}

\proclaim{Proposition 3.4}
Let  $(\Cal X,\zeta )$ be a flat family  of
 polarized K3
surfaces of genus $g$ over the disc $T$, whose central fiber $
(\Cal X_0 , \zeta _0 )$ is a hyperelliptic polarized K3 surface.
Assume furthermore that $\zeta_t$ is a very ample line bundle
on $\Cal X_t$, for all $t \neq 0$. If
\text{$\Cal Y =
\phi _{\zeta} (\Cal X)
\subset {\bold P}_T ^g$ ,}
then the central fiber $\Cal Y_0$ is a K3 carpet. \endproclaim

{\it Proof}.
First, we prove that $\Cal Y_0$ is a carpet.
By Theorem 2.1 we only have to see that through every point $p
\in (\Cal Y_0)_{red} $ there exists a hyperplane $H_p$
such that $\Cal Y_0 \cap H_p$ is a ribbon.
By Bertini's Theorem, we can choose a (generic) hyperplane $H_p$ that
passes through p and whose intersection with $(\Cal Y_0)_{red}$ is a
smooth curve.
$(\Cal Y_0)_{red} \cap H_p$ is a rational normal curve. By Remark 3.1
we know that $\Cal Y$ is flat over $T$. By Lemma 3.2 we know that
$\Cal Y_0
\cap H_p$ is the limit of a family of canonical curves in $\bold P
^{g-1}$, namely, the image of a family of curves whose central
fiber is hyperelliptic (and whose general fiber is not), mapped
by the complete linear series of the relative dualizing sheaf.
Corollary 3.3 tells us that $\Cal Y_0 \cap H_p$
is actually a canonical ribbon.

Second, we prove that the canonical sheaf of $\Cal Y_0$ is
trivial and that the  irregularity of $\Cal Y_0$ is 0. Since $\Cal Y_0$ is
the flat limit of a family of smooth K3 surfaces, $\Cal X
(\Cal O _{\Cal Y_0}) =2$ and therefore $\text h ^2 (\Cal O _{\Cal Y_0})
\geq 1$. Thus, there exists a nonzero global section  $s$
of $\omega _{\Cal Y_0}$. We intend to show that $s$ is
nowhere vanishing. We have the following exact sequence:
$$ 0 \to \omega _{(\Cal Y_0)_{red}} \to \omega _{\Cal Y_0} \to
\Cal I ^* \otimes \omega _{(\Cal Y_0)_{red}} \to 0 $$
that comes from dualizing:
$$ 0 \to \Cal I \to \Cal O _{\Cal Y_0} \to \Cal O _{(\Cal Y_0)_{red}}
\to 0 $$
Therefore H$^0 (\omega _{\Cal Y_0}) = $H$^0(\Cal I ^*
\otimes \omega _{(\Cal Y_0)_{red}})$. This implies that if $s$
vanishes at every closed point of $(\Cal Y_0)_{red}$, then
$s$ is the zero section. Thus, $Z(s) \subsetneqq
(\Cal Y_0)_{red}$. Assume it is not empty and take $p \in
Z(s)$. In the first part of the proof we showed that
the intersection $D $ of the generic
 hyperplane $H_p$ through $p$ with $\Cal Y_0$
 is  a canonical ribbon in $H_p =
\bold P ^{g-1}$.  This  means that $\Cal O _{D}
(1) = \omega _{D} $.  From the adjunction
formula we obtain that $\omega _{\Cal Y_0} |_{D} =
\Cal O _{D}$.
 Hence $s| _{D} \in $H$^0(
\Cal O_{D})$, and since $s$ vanishes at
$p \in
D_{red}$, $s| _{D_{red}}$
must be the zero section,
 but this contradicts the fact that $Z(s) \subsetneqq
(\Cal Y_0)_{red}$. Therefore $Z(s) = \emptyset$ and $\omega
_{\Cal Y_0} \simeq \Cal O_{\Cal Y_0}$. Using  $ \Cal X (\Cal
O_{\Cal Y_0}) = 2$, it follows that h$^1(\Cal O_{\Cal Y_0})
= 0$. \boxit{}

\vskip .5 cm
We will use Proposition 3.4 to prove our main
\proclaim{Theorem 3.5} Any  K3 carpet can be smoothed.
\endproclaim

For the proof of Theorem 3.5 we
need also the following relative version of Proposition 1.5:

\proclaim {Theorem 3.6}
The scheme $U$
  parametrizing (smooth) rational normal \linebreak
scrolls embeds
into the Hilbert scheme of numerical K3 surfaces of degree
$2g-2$ in $\bold P ^g$.
 The image of $U$ by this embedding parametrizes the
K3 carpets in $\bold P ^g$.
In particular, K3 carpets supported on
$S(a-1,b+1)$ lie on the closure of the locus of the Hilbert scheme
parametrizing K3 carpets supported on
$S(a,b)$.
\endproclaim
We use these two propositions to prove Theorem 3.6:

\proclaim{Proposition 3.7} Let $U$ be a smooth variety. If $p:\Cal
S\to U$ is a flat family of rational normal scrolls
inside $\bold P ^g _U$, then there exists a unique family
$\tilde{\Cal S}$ over
$U$, whose fibers are K3 carpets and such that $\tilde{\Cal
S}_{red}=\Cal S$.
\endproclaim

{\it Proof.}
Let $\Cal N$ be the normal bundle of $\Cal S$ inside $\bold P^g_U$
and let $\omega$ denote the relative dualizing sheaf of $\Cal
S/U$, which is in this case a line bundle.
By Proposition 1.5, $p_*(\Cal N \otimes \omega)$ is also a line
bundle.
We claim that $\Cal N \otimes \omega \otimes p^*(p_*(\Cal N \otimes
\omega))^*$ has a nowhere vanishing section $s$ and that $\text H^0(\Cal
N
\otimes \omega \otimes p^*(p_*(\Cal N \otimes
\omega))^*)=s \cdot \text
H^0(\Cal O_U)$. Indeed, by projection formula,
$$\displaylines{\text H^0(\Cal N
\otimes \omega \otimes p^*(p_*(\Cal N \otimes
\omega))^*)
\cr
=\text H^0(p_*(\Cal N \otimes \omega \otimes p^*(p_*(\Cal
N
\otimes
\omega))^*)=
\text H^0(\Cal O_U) \ ,\cr}
$$ and $1 \in \text H
^0(\Cal O_U)$ corresponds to a section $s$ of H$^0(\Cal N
\otimes \omega \otimes p^*(p_*(\Cal N \otimes
\omega))^*)$ which does not vanish identically along  any fiber of
$p$. Since H$^0(\Cal N _{\Cal S_u/\bold P ^g} \otimes \omega
_{\Cal S _u}) = \bold C \cdot
s'$ by Proposition 1.5, where $s'$ is nowhere vanishing section,
it follows that
$s$ is nowhere vanishing. In particular, any nowhere vanishing
section of \text{$\Cal N \otimes \omega \otimes p^*(p_*(\Cal N \otimes
\omega))^*$} is a multiple of $s$ by a global section of $\Cal O_U
^*$. By Lemma 1.3, $s$ defines a double structure $\tilde{\Cal
S}$ on
$\Cal S$ and by the previous observation any other nowhere
vanishing section of
$\Cal N
\otimes
\omega
\otimes p^*(p_*(\Cal N \otimes
\omega))^*$ defines the same double structure.
The ideal sheaf of
$\Cal S$ in $\Cal O _{\tilde {\Cal S}}$ is the line bundle $\omega \otimes
p^*(p_*(\Cal N
\otimes
\omega))^*$. Hence, since $\Cal S$ is flat over $U$, it follows
that $\tilde{\Cal S}$ is also flat over $U$. This implies that
$\tilde{\Cal S}$ is a family of K3 carpets.

Now we prove the uniqueness of $\tilde{\Cal S}$. Let $\tilde{\Cal
S}'$ be a flat family  over $U$, whose
fibers $\tilde{\Cal S}'_u$ are K3 carpets such that $(\tilde{\Cal
S}'_u)_{red}= \Cal S_u$, for all $u \in U$. Using Theorem 2.1
inductively,
(we lift a regular sequence defining the point $u$ in $\Cal O
_{u,U}$ to $\Cal O _{x,\tilde{\Cal S}'}$, where $x$ is any point
in the inverse image of the morphism from $\tilde{\Cal S}'$ to
$U$)
we conclude that $\tilde{\Cal
S}'$ is a double structure on $\Cal S$. This is equivalent to the
data of a vector bundle surjection
$ \Cal N ^* \to \Cal L \to 0$, where $\Cal L$ is a line bundle.
By flatness and because $\tilde{\Cal
S}'$ is a family of K3 carpets, we obtain that $\Cal L |_{\Cal
S_u} = \omega _{\Cal S_u}$ for all $u \in U$. Therefore $(\Cal L
\otimes \omega ^*)|_{\Cal S_u}=\Cal O _{\Cal S_u}$ and $p_*(\Cal
L \otimes \omega ^*)$ is a line bundle.
Moreover, $p^*p_*(\Cal L \otimes \omega^*) \to \Cal L \otimes
\omega ^*$ is a surjective morphism of line bundles
and hence, an isomorphism. Thus $\Cal N \otimes \Cal L = \Cal N
\otimes \omega \otimes p^*p_*(\Cal L \otimes \omega ^*)$. By
hypothesis H$^0(\Cal N \otimes \Cal L)$ contains a nowhere
vanishing section, hence $p_*(\Cal N \otimes \Cal L) = \Cal O
_U$.  By projection formula it follows that $p_*(\Cal L \otimes
\omega ^*) = (p_*(\Cal N
\otimes \omega))^*$  and $\Cal N \otimes \Cal L = \Cal N
\otimes \omega \otimes p^*(p_*(\Cal N
\otimes \omega))^*$. This implies that
$\tilde {\Cal S}' =
\tilde {\Cal S}$.
\boxit{}
\vskip .5 cm

(3.8) \ Recall that  (smooth) rational normal scrolls are
parametrized by a reduced, open subscheme $U$ of the Hilbert
scheme (see, e.g.,
[A]). The subscheme $U$ is stratified as follows (see [A] or [Ha]): the
 scrolls of type $S(a+1,b-1)$, less  {\it balanced}, lie on the
closure of the locus parametrizing scrolls of type $S(a,b)$, more
{\it balanced} (recall that $a \geq b$).

\proclaim{Proposition 3.9} Let $S$ be a rational normal scroll in $\bold
P^N$. The dimension of $\text H^0(\Cal N_{S/\bold P^N})$ is $(N+1)^2-7$ and
$\text H^1(\Cal N_{S/\bold P^N})$ and $\text H^2(\Cal N_{S/\bold P^N})$ vanish.
\endproclaim

{\it Proof.} The statement follows from the exact sequence presenting
$ \Cal N_{S/\bold P^N}$, from the Euler sequence  on $\bold P^N$
and from the sequence relating the tangent bundle of $S$, the
relative tangent bundle of the fibration to $\bold P^1$ and the
pullback of the tangent bundle to
$\bold P^1$. \boxit{}

(3.10)
{\it Proof.}
The scheme $U$ is smooth (by Proposition 3.9 and [S], corollaries
8.5 and 8.6;  see also [A]).  Thus,
we can apply Proposition 3.7 and by the universal property of the
Hilbert scheme, we obtain a morphism $\varphi$ from $U$ to the
Hilbert scheme of numerical K3 surfaces. Let $Z$ be the image of
$\varphi$. The scheme $Z$ parametrizes the K3 carpets inside the
Hilbert scheme. To see that $\varphi$ is an isomorphism onto $Z$
it suffices, since both $U$ and $Z$ are varieties and we are
working over $\bold C$, to show that there exists a morphism
$\Psi$ that is a set-theoretical inverse of $\varphi$. To
construct $\Psi$, consider the pull-back to $Z$ of the universal
family on the Hilbert scheme. The fibers of this pull-back are K3
carpets. If we take the reduced structure of the pull-back,
we  end up with a family of rational normal scrolls over
$Z$. The universal property of the Hilbert scheme gives us the
existence of $\Psi$.

The observation about the stratification of the locus of K3
carpets follows from (3.8).
\boxit{}

\vskip .5 cm
(3.11)
{\it Proof.}
First consider the K3 carpets whose reduced structured
is a rational normal scroll $F$ (embedded in $\bold P ^g$ as a
variety of minimal degree) of type F$_0$,
\dots , F$_4$. The scroll $F$ can be realized as the image of the
morphism induced by the hyperelliptic
linear series of a polarized hyperelliptic K3 surface $(X,
 L)$. We give here a sketch of the construction of
$(X,
 L)$; for more details, see [D] or [R].  Take a curve $C$ in
$|-2K_F|$ with at worst certain mild singularities. Then the
desingularization $X$ of the double cover  of $F$ ramified along
$C$ is a K3 surface. The line bundle $ L $ is the pullback of
$\Cal O_F(1)$. Let $E$ be the
elliptic pencil obtained as pullback of the ruling of $F$. In this
situation the Picard lattice of
$(X,
 L)$ contains a
sublattice generated by $L$ and by
$E$ with intersection
matrix
$$ \pmatrix 2g-2 &2 \cr
2 &0\cr
\endpmatrix$$
Using the fact that the space of periods is a fine moduli
space for polarized, marked K3 surfaces (see [SP]), one can find a
family $(\Cal X , \zeta)$ of polarized K3 surfaces over the analytic
disc $T$, whose central fiber
$(\Cal X_0 ,
\zeta_0)$ is isomorphic to $(X,L)$ and such that $\zeta _t$ is
very ample if $t \neq 0$.  This is achieved  by taking a path in
the period space in this way: the central point corresponds to a
period containing
$E$ and the other points correspond to periods containing
neither $E$ nor any  class with nonpositive intersection with $L$.
Let
$\Cal Y$ be $ \phi _{\zeta} (\Cal X) \subset {\bold P}_T ^g$.
Theorem
1.4 tells us that there exists a K3 carpet structure on
$F$ that can be smooth, namely, $Y_0$. This proves the
theorem in this case, since we know by Theorem 1.3 that
there is a unique K3 carpet structure on any given rational
normal scroll.

We have just proven that K3 carpets on rational normal scrolls of
type $F_0, \dots ,
F_4$ lie on the closure of at least one component parametrizing
smooth K3 surfaces in the Hilbert scheme. By Theorem 3.6, the
remaining K3 carpets lie also in the closure of that (those)
component(s). \boxit{}

\heading {4. The Hilbert scheme near the point of a K3
carpet} \endheading
In this section we study the geometry of the
Hilbert scheme of numerical K3 surfaces (i.e., regular subschemes of
projective space with trivial dualizing sheaf) at the locus
parametrizing K3 carpets. We start by settling  the question of
whether the Hilbert points of the K3 carpets are smooth.

\proclaim{Theorem 4.1}
Let $\tilde S$ be a K3 carpet supported on
$S = S(a,b)$, where $a \geq b$ and $g = a +b$.  The Hilbert
point of
$\tilde S$ is nonsingular iff $0 \leq a-b \leq
2$.
\endproclaim

{\it Proof.}
We have proved in Theorem 3.9 that K3 carpets are smoothable.
 \linebreak
Since the dimension of a component parametrizing smooth
K3 surfaces is \linebreak
dim PGL$(g) + 19 = (g+1)^2 + 18$, a K3
carpet represents a smooth point of the Hilbert scheme iff
$\text h^0(\Cal N_{\tilde S/
\bold P ^g}) = (g+1)^2 + 18$. To compute the
cohomology of
$\Cal N_{\tilde S/
\bold P ^g}$ we tensor the sequence
$$ 0 \to \omega _S \to \Cal O _{\tilde S} \to \Cal O _S \to 0$$
by $\Cal N_{\tilde S/
\bold P ^g}$. Since ${\tilde S}$ is locally a complete
intersection,
the
sheaf $\Cal N_{\tilde S/
\bold P ^g}$ is a vector bundle and
we obtain
$$0 \to\Cal N_{\tilde S/
\bold P ^g} \otimes \omega _S \to \Cal N_{\tilde S/
\bold P ^g} \to \Cal N_{\tilde S/
\bold P ^g} \otimes \Cal O _S \to 0 \ . \eqno (4.1.1) $$
Thus, we have
$$  \chi (\Cal N_{\tilde S/
\bold P ^g})= \chi(\Cal N_{\tilde S/
\bold P ^g} \otimes \omega _S) + \chi (\Cal N_{\tilde S/
\bold P ^g} \otimes \Cal O _S) \ . \leqno (4.1.2)$$
Let $\Cal I _S$ (respectively $\Cal I _{\tilde S}$) be the ideal
of $S$ (respectively $\tilde S$) in $\bold P ^g$. Since $\Cal
I_{\tilde S}/\Cal I_{\tilde S}^2$ is a bundle, taking its dual and
restricting it to
$S$ commute. Hence,
$$ \displaylines{\Cal N_{\tilde S/
\bold P ^g} \otimes \Cal O_S = \Cal Hom _{\tilde S}
(\Cal I_{\tilde S}/\Cal I_{\tilde S}^2, \Cal
O_{\tilde S})
\otimes
\Cal O_S \cr
= \Cal Hom _S(\Cal I_{\tilde S}/\Cal I_{\tilde S}^2 \otimes
\Cal O_S, \Cal O_S)  = \Cal Hom _S(\Cal I_{\tilde S}/\Cal I_{\tilde S}  \Cal
I_S , \Cal O_S).}$$ Therefore $\Cal N_{\tilde S/
\bold P ^g} \otimes \Cal O_S$ sits in the  sequence
$$ 0 \to \Cal Hom _S(\Cal I_{\tilde S}/  \Cal
I_S^2 , \Cal O_S) \to \Cal N_{\tilde S/
\bold P ^g} \otimes \Cal O_S \to \Cal Q \to 0\ . \eqno (4.1.3) $$
{}From sequence
$$ 0 \to \Cal I _{\tilde S} /  \Cal I_S ^2 \to \Cal I _S /\Cal I ^2
_S \to
\Cal I_S/\Cal I_{\tilde S} \to 0 \eqno (4.1.4)$$
we see at once that $\Cal Q$ is a line bundle. We claim that $\Cal
Q =
\omega _S ^{-2}$.  From (4.1.3) it follows that
$\Cal Q =
\bigwedge ^{g-2}  (\Cal N_{\tilde S/
\bold P ^g} \otimes \Cal O_S) \otimes \bigwedge ^{g-3}\Cal I_{\tilde S}/  \Cal
I_S^2 $ . Dualizing  sequence (4.1.4)
and taking wedge we obtain that
$$\bigwedge ^{g-3}\Cal I_{\tilde S}/
\Cal I_S^2 =\omega _S ^* \otimes \bigwedge
^{g-2}\Cal N_{S/
\bold P ^g}^*=\omega_S ^{-2} \otimes \Cal
O_S(-g-1) \ .$$ Using adjunction and  the fact that ${\tilde
S}$ is K3 carpet, it follows that \linebreak
\text{$\bigwedge ^{g-2} \Cal
N_{\tilde S/
\bold P ^g} = \Cal O_{\tilde S}(g+1)$,} and therefore  $\bigwedge ^{g-2} (\Cal
N_{\tilde S/
\bold P ^g} \otimes \Cal O_S) = \Cal O_S(g+1)$, so the claim is
clear. Therefore we obtain the following exact
sequences:
$$ \matrix 0&\to &
\Cal Hom _S(\Cal I_{\tilde S}/  \Cal
I_S^2 , \Cal O_S) &\to &\Cal N_{\tilde S/
\bold P ^g}\otimes \Cal O _S &\to& \omega_S ^{-2}&
\to& 0\cr
&&&&&&&(4.1.5)
\cr 0 &\to &\Cal Hom _S(\Cal I_{\tilde S}/  \Cal
I_S^2 , \Cal O_S) \otimes \omega _S&\to&
\Cal N_{\tilde
S/
\bold P ^g}\otimes \omega _S &\to& \omega_S ^*
&\to &0\cr
&&&&&&&(4.1.6)
 \endmatrix$$
and from (4.1.4) we obtain
$$ \matrix
  0 &\to &\omega_S ^*&\to &\Cal N_{ S/
\bold P ^g} &\to &\Cal Hom _S(\Cal I_{\tilde S}/  \Cal
I_S^2 , \Cal O_S) &\to & 0 \cr
&&&&&&& (4.1.7)\cr 0 &\to& \Cal O_S &\to
&\Cal N_{ S/
\bold P ^g}\otimes \omega _S &\to& \Cal Hom _S(\Cal I_{\tilde S}/
\Cal I_S^2 , \Cal O_S) \otimes \omega _S &\to &0\cr
&&&&&&& (4.1.8)
 \endmatrix$$
Using (4.1.8) and  Proposition 1.5, it follows that $ \text H
^1(\Cal Hom _S(\Cal I_{\tilde S}/
\Cal I_S^2 , \Cal O_S) \otimes \omega _S) = 0$. From
(4.1.7) and Proposition 3.7, and  the fact that  $\text H
^2(\omega_S ^*)=0$, it follows that $\text H ^1(\Cal Hom
_S(\Cal I_{\tilde S}/  \Cal I_S^2 , \Cal O_S) )$ vanishes.
Therefore
$$ \displaylines {\text h ^0(\Cal N_{\tilde S/
\bold P ^g}\otimes \Cal O _S) = \text h ^0(\Cal Hom _S(\Cal I_{\tilde S}/
\Cal I_S^2 , \Cal O_S) ) + \text h ^0(\omega_S ^{-2})
\cr = \text h
^0(\Cal N_{ S/
\bold P ^g}) - \text h
^0(\omega_S ^*) + \text h ^1(\omega_S ^*) + \text h ^0(\omega_S
^{-2})}$$
and
$$\displaylines {\text h ^0(\Cal N_{\tilde
S,
\bold P ^g}\otimes \omega _S) =  \text h ^0(\Cal Hom _S(\Cal I_{\tilde S}/
\Cal I_S^2 , \Cal O_S) \otimes \omega _S) + \text h ^0(\omega_S
^{*}) \cr = \text h ^0(\Cal N_{ S/
\bold P ^g}\otimes \omega _S) - 1 + \text h ^0(\omega_S ^{*}) =
\text h ^0(\omega_S ^{*}),}
$$ by Proposition 1.5. By Proposition 3.7, $\text H ^2(\Cal
N_{\tilde S/
\bold P ^g}) = 0$. From this and (4.1.7) it follows that $\text H ^2(\Cal Hom
_S(\Cal I_{\tilde S}/  \Cal I_S^2 , \Cal O_S) ) $ vanishes.
Using that $\text H ^2(\omega_S ^{-2})$ \linebreak
and $\text H ^1(\Cal Hom _S(\Cal I_{\tilde S}/
\Cal I_S^2 , \Cal O_S) )$
vanish, from (4.1.5)  we obtain that $ \text H ^1(\Cal N_{\tilde S/
\bold P ^g}\otimes \Cal O _S) = \text H ^1(\omega_S ^{-2})$ and
that $\text H ^2(\Cal N_{\tilde S/
\bold P ^g}\otimes \Cal O _S) = 0$. Analogously, from Proposition
1.5 and  sequence (4.1.8) it follows that
\text{$\text H ^2(\Cal Hom _S(\Cal I_{\tilde S}/  \Cal I_S^2 , \Cal
O_S)\otimes\omega _S)$} vanishes. Using  that $\text H ^2(\omega_S ^{*})$
and \text{$
\text H ^1(\Cal Hom _S(\Cal I_{\tilde S}/  \Cal I_S^2 , \Cal O_S)
\otimes \omega _S)$}
vanish, from (4.1.6) we obtain that \linebreak
$\text H^1(\Cal N_{\tilde S/
\bold P ^g} \otimes \omega _S) = \text H ^1(\omega_S ^{*})$ and
that $\text H ^2(\Cal N_{\tilde S/
\bold P ^g} \otimes \omega _S) = 0$. Therefore we can rewrite (4.1.2) as:
$$\displaylines{\text h ^0(\Cal N_{\tilde S/
\bold P ^g}) - \text h ^1(\Cal N_{\tilde S/
\bold P ^g}) = \text h ^0(\Cal N_{ S/
\bold P ^g}) + \text h ^0(\omega_S ^{-2}) - \text h ^1(\omega_S
^{-2}) \cr
= \text h ^0(\Cal N_{ S/
\bold P ^g}) + \chi (\omega_S ^{-2}).\cr} $$
In Proposition 3.7 we show that the dimension of $\text H ^0(\Cal
N_{ S/
\bold P ^g})$ is $(g+1)^2-7$. By Riemann-Roch one
obtains that $\chi (\omega_S ^{-2}) = 25$. Thus, $\text h ^0(\Cal N_{\tilde S/
\bold P ^g}) - \text h ^1(\Cal N_{\tilde S/
\bold P ^g}) = (g+1)^2+18$ and from this it follows that $\tilde S$
represents a nonsingular point of the Hilbert scheme iff $\text h ^1(\Cal
N_{\tilde S/
\bold P ^g}) = 0$. From sequence (4.1.1) we get
$$
\text H^1 (\Cal N_{\tilde S/
\bold P ^g} \otimes \omega _S)  \to \text H^1(\Cal N_{\tilde S/
\bold P ^g}) \to \text H^1(\Cal N_{\tilde S/
\bold P ^g} \otimes \Cal O _S) \to 0 \ ,$$
hence the key point is to compute the dimension of $
\text H ^1(\Cal N_{\tilde S/
\bold P ^g}\otimes \Cal O _S) = \text H ^1(\omega_S ^{-2})$ and of
$\text H^1(\Cal N_{\tilde S/
\bold P ^g} \otimes \omega _S) = \text H ^1(\omega_S ^{*})$.
Pushing down to $\bold P^1$, we obtain that
$$\displaylines{\hfill \text H ^1(\omega_S
^{*}) = \text H^1(\Cal O _{\bold P^1} (a-b +2) \oplus \Cal O _{\bold
P^1} ( 2) \oplus \Cal O _{\bold P^1} (b-a +2)) \hfill \cr
\rlap{\quad
\text{and that}}\hfill
\text H ^1(\omega_S ^{-2})=
\text H^1(\bigoplus_{i=-2}^2 \Cal O _{\bold P^1} (i(b-a) +4))\
.\hfill\cr}$$ Therefore, if
$0 \leq a-b \leq 2$, both $\text H ^1(\omega_S ^{*})$ and $\text H
^1(\omega_S ^{-2})$ are zero. Hence \text{$\text H^1(\Cal N_{\tilde S/
\bold P ^g})$} vanishes and $\tilde S$ corresponds to
a nonsingular point of the Hilbert scheme. On the other hand, if
$a-b > 2$, the group $\text H
^1(\omega_S ^{-2})$  does not vanish and neither does
$\text H^1(\Cal N_{\tilde S/
\bold P ^g})$.
\quad \boxit{}

\vskip .5 cm
As consequence of Theorem 4.1 we know that K3 carpets on rational
normal scrolls of type
$F_0, F_1, F_2$ belong only to one component of the Hilbert
scheme of numerical K3 surfaces, and by Theorem 3.9, we know that
the general point of that component is a smooth K3 surface. By
using the smoothing constructed in the proof of 3.9 we are able to
identify the component in question. The same construction allows
us to prove that a K3 carpet contained in $\bold P ^g$, when $g
\equiv 1
\ (4)$, and with reduced part isomorphic to the ruled surface
$F_4$, belongs to two components of the Hilbert scheme. This fact
provides a geometric explanation for the nonsmoothness of its
Hilbert point.

\proclaim{Theorem 4.2} The K3 carpets supported on rational
normal scrolls of type $F_0, F_1$ (and therefore any K3 carpet)
belong to the ``prime" component of the Hilbert scheme of
numerical K3 surfaces.
\endproclaim

{\it Proof.} Let $X$ be a hyperelliptic K3 surface mapping
generically
$2:1$ to
$F_0$ or $F_1$. If $X$ maps to $F_0$, the
Picard group of   $X$  contains a sublattice
generated by two elliptic pencils
$E_1$ and
$E_2$.  This sublattice has intersection matrix
$$ \pmatrix
0 & 2 \\
2 & 0 \\
\endpmatrix $$
If $X$ maps to $F_1$, the
Picard group of   $X$  contains a sublattice generated by an elliptic pencil
$E$ and by a rational
{\it nodal} curve $R$. This sublattice has intersection matrix

$$ \pmatrix
0 & 2 \\
2 & -2 \\
\endpmatrix $$
 It is easy to check
that these sublattices are primitive and, in particular, that
$L_n = E_i + nE_j$ is primitive for all $n \geq 1$
and that
$L_n = R + nE$ is primitive for all $n \geq 2$ . The line
bundles $L_n$ are the hyperelliptic
line bundles which give a generically $2:1$ map from $X$ to a
rational normal scroll of type $F_0$ or $F_1$. Using the same
reasoning as in the proof of 3.7 we can construct a family $(\Cal X,
\zeta)$ of polarized K3 surfaces whose central fiber is
isomorphic to $(X,L_n)$ and whose general fiber $(\Cal X _t,
\zeta_t)$  is a K3 surface such that $\text{Pic}(\Cal X
_t)$ is generated by $\zeta_t$. Therefore we can
construct a smoothing of the K3 carpet supported
on a rational normal scroll of type $F_0$ or
$F_1$ such that the Picard group of the general fiber is generated
by the hyperplane class.
\boxit{}

\proclaim{Theorem 4.3} Let $g$ be  greater than $9$
and congruent to $1$ modulo $4$. The  K3 carpets inside
$\bold P ^g$, supported on a
rational normal scroll $S$ of type
$F_4$  belong to two components of the Hilbert scheme. One of
them is the``prime" component. A general point of the other
component corresponds to a smooth K3  surface with Picard
number one but with hyperplane class divisible by two.
\endproclaim

{\it Proof.}
The Picard group of a hyperelliptic K3 surface $X$ mapping
generically
$2:1$ to
$F_4$  has a sublattice
generated by an
elliptic pencil
$E$ and by a rational
 nodal curve $R$.  This sublattice has intersection matrix
$$ \pmatrix
0 & 1 \\
1 & -2 \\
\endpmatrix $$
 (see [D] for details). The hyperelliptic  line bundles
mapping $X$ generically $(2:1)$ to a rational normal scroll of
type $F_4$ are the line bundles $L_n = 2R + nE$ for all $n \geq
5$. If $n$ is even, the line bundle $L_n$ is not primitive, but
the double of other line bundle. Therefore we can construct in
that case a smoothing of
$\tilde S$ with the following property: the general fiber has
Picard number one but its hyperplane class does not generate the
Picard group, but it is divisible by two. Thus the general fiber
does not belong to the prime component. The hypothesis on $g$
being congruent to $1$ modulo $4$ comes in at this point, because
in that case $n$ is even ($n=\frac {g+3} {2}$).
\boxit{}

\vskip .5 cm

We will devote the rest of the section to describe the
deformation of K3 carpets to the union of two scrolls. Recall that
the union of two rational normal scrolls of dimension 2 along a
(reduced, but maybe reducible) elliptic curve, anticanonical with
respect to both of them, has the numerical invariants of a K3
surface. Ciliberto, Lopez and Miranda prove in [CLM] that
those unions of scrolls having  smooth double locus (note that this
condition forces the reducible K3 to be  a union of two copies of
$F_0$, $F_1$ or $F_2$) are smoothable. In fact, since
any union of two rational normal scrolls along a reduced anticanonical
curve  can be deformed to a union of two scrolls with
smooth double locus, it follows that any union of two scrolls along
a reducible and irreducible anticanonical elliptic curve is
smoothable. Thus, which follows provides another, more indirect,
proof of the smoothing of K3 carpets.

\vskip .1 cm
\proclaim{Theorem 4.4}
The locus of K3 carpets lies on the closure of the locus
parametrizing unions of two scrolls. Both loci lie
on the closure of the open subscheme parametrizing smooth K3
surfaces in the prime component.
\endproclaim

{\it Proof.}
Let $S$ be  a rational normal scroll. Let $C$ be a curve in the
linear equivalent class of the anticanonical divisor. The curve $C$
induces an embedding
$$ 0 \to \Cal N _{S,\bold P^g} \otimes \omega \to \Cal N _{S,\bold
P^g}
$$ and the image of the generator of H$^0(\Cal N _{S,\bold P^g}
\otimes \omega)$ in H$^0(\Cal N _{S,\bold P^g}
)$ corresponds to a first order deformation of $S$ in $\bold P
^g$, keeping $C$ fixed. Since h$^1(\Cal N _{S,\bold P^g}
\otimes \omega) = 0$ by Proposition 1.5, this first order
deformation extends to a deformation of $S$ over a smooth affine
curve
$U$, keeping $C$ fixed. We will call this deformation $\Cal S
_1$ and by an abuse of notation, we will denote its central fiber
by $S$. Consider now another deformation $\Cal S_2$ fixing $C$
(e.g., the trivial deformation
$S
\times U \subset \bold P ^g_U$). The family $\Cal S_1 \cup \Cal S_2$
is flat over $U$ and the general fiber is the union of two scrolls.
We claim that the central fiber is a K3 carpet. Note that $\Cal
S_1 \cap \Cal S_2 = S \cup (C \times T)$. For any point $x
\in S$ we choose a hyperplane $H_1$ passing through $x$ such
that $D := S \cap H$ is a smooth rational normal curve and such
that $(H_1 \times U) \cap (\Cal S_1 \cap \Cal S_2)$ is induced locally by a
non-zero-divisor on
$\Cal O_{\Cal S_1 \cap \Cal S_2}$. The scheme  $(\Cal S_1 \cup \Cal
S_2) \cap (H_1 \times U)$ is flat and equal to $(\Cal S_1 \cap (H_1
\times U)) \cup (\Cal S_2 \cap (H_1 \times U))$. Denote $\Cal S_i
\cap (H_1
\times U)$ by $\Cal S'_i$. Now $\Cal
S'_1 \cap \Cal S'_2 = D \cup ((C \cap H_1) \times T)$ and
through any point $y \in S \cap H_1$, we choose a hyperplane
$H_2$ such that
$D
\cap H_2$ consists of  distinct points and $(H_2 \times U) \cap (\Cal S_1'
\cap \Cal S_2')$ is induced by a non-zero-divisor of $\Cal O_{\Cal S'_1 \cap
\Cal
S'_2}$. Again, the family $(\Cal S'_1 \cup \Cal
S'_2) \cap (H_2 \times U)$ is flat and equal to $(\Cal S'_1 \cap
(H_2
\times U)) \cup (\Cal S'_2 \cap (H_2 \times U))$. The general
fiber of $(\Cal S'_1 \cup \Cal
S'_2) \cap (H_2 \times U)$ consists of $2g-2$ distinct points and
the central fiber is supported on $g-1$ distinct points. Now the
proof follows the same path as the proof of Proposition 3.4. By the
same degree considerations, the central fiber of  $(\Cal S'_1 \cup \Cal
S'_2) \cap (H_2 \times U)$ is a 0-dimensional ribbon. The central
fiber of $\Cal S'_1 \cup \Cal S'_2$ is also a ribbon by Theorem
2.1 and Lemma 3.2. In fact, it is a canonical ribbon, because it
is a nondegenerate ribbon of degree $2g-2$ in $\bold P ^{g-1}$.
Again by Theorem 2.1 and Lemma 3.2, we obtain that the central
fiber of $\Cal S_1 \cup \Cal S_2$ is a carpet, and adjunction
implies that it is a K3 carpet.
\quad \boxit{}
\vskip .5 cm

(4.5) An example of this degeneration can be constructed
explicitly in the following way: let $S$ be  a rational normal
scroll in $\bold P ^g$, $g \geq 4$. Let $C_0$ be
the minimal section of $S$. Fix a smooth section $C'$ not
intersecting $C_0$. Let
$\phi$ be the morphism from $C_0$ to
$C'$ defined by the fibers of $S$. Fix three points $a_0, b_0, c_0$ in $C_0$
and let $\phi (a_0) = a', \phi (b_0) = b'$ and $ \phi (c_0) = c'$.
Define $\phi _x : C_0 \to C'$ as the morphism that sends  $a_0$ to
$a'$, $b_0$ to $b'$ and $c_0$ to $x$. Let $(D,d)$ be a smooth
projective curve and
$f: (D,d)
\to (C',c')$  a covering of $(C',c')$. Let $\Psi _f$ be
defined as follows:
$$\matrix
\Psi_f : & C_0
\times D &\to & C'\\
& (t,y) & \mapsto &\phi _{f(y)}(t) \\
\endmatrix$$
The morphism $\Psi _f$ defines  a family $S_f$ of rational
normal scrolls, parametrized by $D$.
Each member of the family contains the
reducible elliptic curve \linebreak
\text{$C := C_0 \cup C'
\cup
\overline
{a_0 a'} \cup \overline {b_0 b'}$,}  which is
anticanonical in each scroll. Choosing $f_1, f_2
: (D, d) \to (C',c')$, $f_1 \neq f_2$ we obtain
$\Cal S_1 = S_{f_1}$ and $\Cal S_2 = S_{f_2}$ and
a family $\Cal S_1 \cup \Cal S_2$ like in the proof
of 4.4.  In fact, all this construction takes
place inside the join $\Sigma$ of
$C_0$ and
$C'$. We will study the situation in more detail:

\proclaim{Proposition 4.6} Let $S$ and $\Sigma$ be as in (4.5). The
3-fold $\Sigma$ is Fano (i.e., some power of its anticanonical divisor is
ample). The K3 carpet supported on
$S$ is a member of
$|-\omega _\Sigma|$. All members of $|-\omega _\Sigma|$ are singular.
\endproclaim

{\it Proof.}
First we show that $\Sigma$ is Fano. Let $\pi:\Gamma \to \Sigma$ be
the blowing up of
$\Sigma$ along
$C_0$ and $C'$. Let
$n_0$ and
$n'$ be the degrees of the two rational normal curves and let $\Cal
E$ be equal \linebreak
to \text{$
\Cal O _{\bold P^1
\times
\bold P ^1} (n_0,0)
\oplus \Cal O _{\bold P^1 \times \bold P ^1} (0,n')$.} Note that
$\Gamma$ is
$\bold P (\Cal E)$. Let $p$ denote the map from $\Gamma$ to $\bold
P^1 \times \bold P ^1$. Note that if both
$n_0$ and
$n'$ are greater than
$1$, $\Gamma$ is a minimal desingularization of $\Sigma$. In any
case, $\Gamma$ is smooth. It is easy to compute the Chow ring of
$\Gamma$: the generators of A$^3(\Gamma)$ are the pull back of the
first ruling of $\bold P^1 \times \bold P ^1$, that we will denote
by $A$, the pullback of the second ruling, denoted by $B$, and $H$,
the divisor corresponding to $\Cal O _{\bold P (\Cal E
)} (1)$. The class $A$ is represented by a
ruled surface of type $F_{n_0}$ and $B$ is represented by  a ruled
surface of type $F_{n'}$. The generators of A$^2(\Gamma)$ are the
class $C_1$ of the minimal section of $A$, the class $C_2$ of the
minimal section of $B$, and the class $f$ of  a fiber of $p$.
The intersection in the Chow ring is given by the following
matrices (bases of A$^3(\Gamma)$ and A$^2(\Gamma)$ are ordered as introduced
before):
$$\pmatrix n'C_1 + n_0C_2 + n_0n'f & C_1 + n_0f & C_2 + n'f \cr
C_1 + n_0f & 0 & f \cr
C_2 + n' f & f & 0 \cr
\endpmatrix$$
and
$$\pmatrix
0 & 0 & 1 \cr
0 & 1 & 0 \cr
1 & 0 & 0 \cr
\endpmatrix$$
Using this information and adjunction we compute the class of the
canonical divisor of $\Gamma$:
$$ K_\Gamma \equiv (n' -2) A + (n_0 - 2) B - 2H$$
In fact, the previous equality is up to linear equivalence, since
$\text H ^1(\Cal O _\Gamma) = 0$. This implies that $\Gamma$ is a Fano
3-fold if both $n_0$ and $n'$ are less or equal than $2$.
Otherwise, the anticanonical divisor has negative intersection with
$C_1$, with $C_2$, or with both of them. But $C_1$ and $C_2$  are contracted
by the morphism from $\Gamma$ to $\Sigma$ so in any case $\Sigma $ is a Fano
3-fold.

For the second claim, let $\tilde S$ be a carpet inside $\Sigma$,
supported on a rational normal scroll $S$. Let $\tilde S'$ be the
strict transform of $\tilde S$ in $\Gamma$ and let $E_1$ (respectively
$E_2$) be the exceptional divisor corresponding to the rational
normal curve of degree $n_0$ (respectively of degree $n'$). We
would like to prove that
$\pi ^*\Cal O (K_\Sigma + \tilde S) = \Cal O_\Gamma$, because that would
imply that $\tilde S \sim -K_\Sigma$ ($\pi_*\Cal O_\Gamma = \Cal O_{\Sigma}$
because
$\Sigma$ is normal). The 3-fold
$\Sigma$ is, locally along the
rational normal curves $C_0$ and $C'$, formally isomorphic to the product of a
line with the blowing down of a surface at a $-n_0$ or $-n'$ -curve.
check!  Therefore, $\Sigma$ is $\bold Q$-Gorenstein and
$$ \displaylines{\rlap{(4.6.1)}\hfill\pi ^* K_\Sigma \equiv K_\Gamma + \frac
{n' - 2} {n'} E_1 +
\frac {n_0 - 2} {n_0} E_2\hfill\cr
\hfill\equiv (n' -2) A + (n_0-2) B - 2H +
\frac {n' - 2} {n'} E_1 + \frac {n_0 - 2} {n_0} E_2 \ .\hfill}$$ On
the other hand, $\tilde S '$ is a carpet inside $\Gamma$ supported
on a rational normal scroll. The class of the rational normal scroll
is $A+B$. Hence the class of $\tilde S'$ is $2A+2B$. Thus we can
write
$$\pi ^*\tilde S \equiv 2A + 2B + \frac {2} {n'}
E_1 +
\frac {2} {n_0} E_2 \ .\leqno(4.6.2) $$
Putting (4.6.1) and (4.6.2) together we
see at once that
$$\pi
^*(K_\Sigma + \tilde S) \sim n' A + n_0B -2H + E_1 - E_2 \ .$$
To finish
the computation we write $E_1$ and $E_2$ in terms of $A$,$B$ and
$H$. For that we compute the intersections of $E_1$ and $E_2$
with $C_1$, $C_2$ and $f$:
$$ \matrix E_1 \cdot C_1 &\equiv& (E_1 |_ A) |_{C_1} &
\equiv &(C_1 +
n_0f) \cdot C_1
&\equiv &0 \cr
E_1 \cdot C_2 &\equiv &(E_1 |_ B) |_{C_2} &\equiv& C_2 ^2
&\equiv &-n' \cr
E_1 \cdot f &\equiv& 1\cr
\endmatrix$$
and analogously $E_2 \cdot C_1 \equiv -n_0$, $E_2 \cdot C_2 \equiv 0$ and
$E_2 \cdot f \equiv 1$. From this it follows that $E_1 \sim -n'A +
H$ and $E_2 \sim -n_0B + H$.  Therefore $\pi ^*\Cal O (K_\Sigma
+\tilde S) = \Cal O_\Gamma$ as we wanted.

To prove that all elements of $|-K_\Sigma|$ are singular, first we
note that no element of $|-K_\Sigma|$ is contained in the smooth locus
of $\Sigma$. This follows from the fact that $-K_\Sigma$ is
ample, and hence must intersect positively any irreducible
curve, in particular those in the singular locus. To conclude our
argument, we consider two cases, the case when $-K_\Sigma$ is
Cartier and the case when it is not. If $-K_\Sigma$ is Cartier, all
members of  $|-K_\Sigma|$ are singular because they have nonempty
intersection with the singular locus of
$\Sigma$. If $-K_\Sigma$ is not Cartier, let us assume that there
exist a smooth member
$R$ of $|-K_\Sigma|$. The intersection of $R$ and the singular locus of
$\Sigma$ cannot
have dimension $0$. Indeed. Assume the contrary. Let
$R'$ be the strict transform of $R$ by $\pi$. Away from  $C_0$
and $C'$ the sheaf
$\pi_*
\Cal O_\Gamma(R')$ is a line bundle. Using the theorem on formal functions
one sees that $\pi_*
\Cal O_\Gamma(R')$ is in fact a line bundle everywhere (see [K] for a
similar situation in dimension 2). Then, since $\pi_*
\Cal O_\Gamma(R') = \Cal O_\Sigma(R) = \omega_\Sigma^*$ outside from a locus of
codimension 2, it follows that $\pi_*
\Cal O_\Gamma(R')  = \omega_\Sigma^*$ everywhere and hence $\omega_\Sigma^*$ is
a
line bundle, which is a contradiction. Therefore $R$ contains at least one
component of the singular locus (one of the two rational normal curves $C_0$
and $C'$). Since $R$ is smooth it follows that
$$0 \sim \pi ^*(K_\Sigma + R) \equiv K_\Gamma + R' + \frac {n'-2+c} {n'} E_1 +
\frac {n_0 - 2 +d} {n_0} E_2$$
for some $c , d = 0,1$. Since $R$ is smooth along $C_0$ and $C'$, $R$ and
$R'$ are isomorphic. On the other hand $R$ is a flat deformation of the
K3 carpet $\tilde S$, hence, by adjunction $ (n'-2+c)\cdot n_0  (E_1
\cdot R') + (n_0-2+d)\cdot n' (E_2 \cdot R')\equiv 0$ in
$\text{Num}(\tilde R)$, which is not possible because $\text
H^1(R')=0$.

A quicker way to see that all elements of $|-K_\Sigma|$ are
singular is by using the fact that $\text
H^0(\omega_\Sigma^*) = \text H^0(\omega_{C_0}^*) \otimes \text
H^0(\omega_{C'}^*)$ (c.f. [BE], \S 8). Hence an element of
$|-K_{\Sigma}|$ is the union of two cones over $C_0$ and two cones
over $C'$ and therefore is singular along both $C_0$ and $C'$.
\boxit{}
\vskip .35 cm

The same  argument proves that the union
of two scrolls, both of which are obtained by joining
corresponding points of $C'$ and $C_0$, is in the class of the
anticanonical divisor of $\Sigma$
and therefore the
deformation constructed in (4.5) comes from deforming the K3 carpet inside
 $|-K_\Sigma|$. The Proposition tells us also that it is not
possible to construct a smoothing of the K3 carpet inside
$|-K_\Sigma|$.
\vskip .35 cm
(4.7) When $a-b \leq 2$, this construction fits into a more
general one: consider a smooth elliptic normal curve $E$ in $\bold
P ^g$. Let
$\Omega$ be its 2-secant variety. The variety $\Omega$ is
a ``fake"
Calabi-Yau 3-fold (its dualizing sheaf is trivial and the
intermediate cohomology of its structure sheaf vanishes, but its
desingularization is a projective bundle over $S_2(E)$, hence it
has negative Kodaira dimension) see [GP] for details. It is singular along
$E$. A
$g^1_2$ on
$E$ defines a rational normal scroll containing $E$ as a member of the
anticanonical class. If we consider two families of $g^1_2s$
specializing to a given one (the one defining the scroll on which
our K3 carpet is supported) we obtain again a family like in the
proof of 4.5. To go from this picture to the previous one, we
just degenerate $E$ to $C$. The 3-fold $\Omega$ degenerates
to a reducible variety, one of whose components is $\Sigma$. Finally, we can
identify the degenerations of the
$g^1_2s$ as pencils having degree $1$ on $C_0$ and $C'$ and degree
$0$ on
$\overline {a_0 a'}$ and $\overline {b_0 b'}$.
For example, in $\bold P ^4$, the variety $\Omega$ is a quintic 3-fold (in
this case, since $\Omega$ is a hypersurface, one can easily check that is
Calabi-Yau). The degeneration of $\Omega$ consists  of $\Sigma$,
which is a quadric cone, and three hyperplanes.

\heading  {\bf REFERENCES} \endheading
\roster
\item"[ACGH]" E. Arbarello, M. Cornalba, P.A. Griffiths, J.
Harris, {\it Geometry of Algebraic Curves, Vol. I}, Springer, New
York 1985.

\item"[A]" P. Azcue, {\it On the dimension of the
Chow varieties}, Harvard Thesis (1992).

\item "[BE]" D. Bayer \& D. Eisenbud,   {\it
Ribbons and their canonical embeddings}, Trans.
Amer. Math. Soc. {\bf 347} (1995),
719-756.

\item"[CLM]" C. Ciliberto, A. Lopez \& R. Miranda,
{\it Projective degenerations of K3
surfaces, Gaussians maps, and Fano threefolds},
Invent. math. {\bf 114} (1993),
641-667.

\item"[D]" I. V. Dolgachev, {\it On special
algebraic K3 surfaces. I}, Math. USSR Izvestija
{\bf 7} (1973), 833 - 846.

\item"[E]" D. Eisenbud, {\it Green's conjecture; an orientation
for algebraists}, Sundance 91: Proceedings of a Conference on Free
resolutions in Commutative Algebra and Algebraic Geometry, Jones
and Barlett, 1992, 51-78.

\item"[EG]" D. Eisenbud \& M. L. Green, {\it Clifford indices of
ribbons}

\item"[F]" L.-Y. Fong, {\it Rational ribbons and
deformations of hyperelliptic curves},
J. Algebraic Geom. {\bf 2} (1993), 295-307.

\item"[GP]" M.Gross \& S.Popescu, {\it Equations of $(1,d)$ polarized abelian
surfaces}, preprint.

\item"[Ha]" J. Harris, {\it A bound on the
geometric genus of projective varieties}, Ann.
Sci.
Norm. Sup. Pisa {\bf IV-VIII} (1981), 35-68.

\item"[H]" R. Hartshorne {\it Algebraic Geometry}, Springer,
Berlin, 1977.

\item"[HV]" K. Hulek \& A. Van de Ven, {\it The
Horrocks-Mumford bundle and the Ferrand
construction}, Manuscripta math. {\bf 50} (1985),
313-335.

\item"[K]" J. Koll\'ar, {\it Contractibility and minimal model program for
log surfaces},\linebreak preprint (1995).

\item"[R]" M. Reid, {\it Hyperelliptic linear
systems on a K3 surface}, J. London Math.
Soc. {\bf (2), 13} (1976), 427-437.

\item"[SP]" S\'eminaire Palaiseau, {\it
G\'eom\'etrie des surfaces K3:  modules et
p\'eriodes}, \linebreak
Ast\'erisque
{\bf 126}, Soc. Math. France
(1985).

\item"[S]" E. Sernesi, {\it Topics on families of projective schemes},
Queen's Papers in Pure and Appl. Math. 73, Queen's Univ., Kingston Canada,
1986.

\endroster

\enddocument